\documentclass[12pt,eqright]{amsart}
\usepackage{mathrsfs}
\usepackage{amsmath,amssymb}
\usepackage[displaymath,pagewise]{lineno}
\usepackage{bm}
\usepackage{amsfonts}
\usepackage{epsf}
\usepackage{graphicx}
\newtheorem{theorem}{Theorem}[section]
\newtheorem{lemma}[theorem]{Lemma}

\newtheorem{remark}[theorem]{Remark}

\newtheorem{algorithm}{Algorithm}[section]
 \theoremstyle{definition}

\numberwithin{equation}{section}

\newcommand{\bc}{\begin{center}}
\newcommand{\ec}{\end{center}}
\newcommand{\be}{\begin{eqnarray}}
\newcommand{\ee}{\end{eqnarray}}

\newcommand{\ben}{\begin{eqnarray*}}
\newcommand{\een}{\end{eqnarray*}}

\newcommand{\Vvert}{|\hskip-0.8pt|\hskip-0.8pt|}
\newcommand{\Om}{\Omega}
\newcommand{\al}{\alpha}

\newcommand{\pa}{\partial}

\newcommand{\na}{\nabla}

\def\x{\times}

\def\na{\nabla}

\def\ds{\,ds}

\def\cE{\mathcal{E}}

\def\cT{\mathcal{T}}
\def\R{\mathbb{R}}

\def\bpsi{\boldsymbol{\psi}}

\def\div{\operatorname{div}}

 \DeclareMathOperator{\res}{R}
\DeclareMathOperator{\id}{id}

\DeclareMathOperator{\card}{card}
\DeclareMathOperator{\osc}{osc}

\DeclareMathOperator{\consis}{\kappa}

\title[Convergence and optimality of the ANFEM]
{Convergence and optimality of  the \\
  adaptive nonconforming linear element method \\for the Stokes
  problem} \author[J.~Hu]{Jun Hu} \address{LMAM and School of
  Mathematical Sciences, Peking University, Beijing 100871,
  P. R. China} \email{hujun@math.pku.edu.cn} \author[J.~Xu]{Jinchao
  Xu} \address{The School of Mathematical Sciences and  the Beijing International Center for Mathematical Research,  Peking
  University, and Department of Mathematics, Pennsylvania State
  University, University Park, PA 16801} \email{xu@math.psu.edu}

\thanks{The first author was supported by NSFC  10971005,
  and in part by  NSFC
  11031006.  The  second author was supported in part, by NSFC-10528102,
NSF DMS 0915153,  and DMS 0749202,  and by the PSU-PKU Joint Center for
Computational Mathematics and Applications}

\date{\today}
\begin{document}
\begin{abstract}
  In this paper, we analyze the convergence and
  optimality of a standard adaptive nonconforming linear element
  method for the Stokes problem.  After establishing a special
  quasi--orthogonality property for both the velocity and the pressure
  in this saddle point problem, we  introduce a new prolongation
  operator to carry through the discrete reliability analysis for the
  error estimator. We  then use a specially defined interpolation
  operator to prove that,  up to oscillation,  the  error can
  be bounded by the approximation error within a properly defined
  nonlinear approximate class.  Finally, by introducing a new
  parameter-dependent error estimator, we  prove the
  convergence and optimality estimates.

\end{abstract}
\maketitle
\section{Introduction}
The adaptive finite element method plays an important role in the
numerical solution for partial differential equations
\cite{AinsworthOden00,Verfurth96}.  The convergence and optimality of
the adaptive method have been much studied in recent years.  For the
Poisson equation and its variants, the theory is well--developed
\cite{CarstensenHoppe05b,CasconKreuzerNochettoSiebert07,ChenHolstXu07,Dolfler96,MekchayNochetto05, MorinNochettoSiebert00,
  Stevenson06,Stevenson05}.
However, for many other important
problems this is not the case. Among these under studied problems is the Stokes problem, the main subject of this  paper.

The convergence analysis of the adaptive finite element method of the Poisson equation  is based on the orthogonality property
\cite{CasconKreuzerNochettoSiebert07,Dolfler96,
  MekchayNochetto05,MorinNochettoSiebert00}, such orthogonality  can be weakened
  to some quasi--orthogonality for the nonconforming and
mixed methods \cite{BeckerMao2008,BeckerMaoShi2010,CarstensenHoppe05a,
  CarstensenHoppe05b,CarstensenRabus2010,ChenXuHoppe2010,ChenHolstXu07,HSX10,HuangHuangXu2010,MaoZhaoShi2010,Rabaus2010}.
The Stokes problem, as a saddle point problem with two variables
(velocity and pressure),  lacks the usual orthogonality or
quasi--orthogonality that holds for the positive and definite problem.  As a
result, it is not obvious how the technique for nonconforming and
mixed methods for the Poisson equation can be carried over to the
Stokes problem.  Although the mixed formulation of the Poisson
equation is also a saddle point problem,  analyses of this formulation's
convergence and optimality \cite{BeckerMao2008,CarstensenRabus2010,ChenHolstXu07} are not so
different from that for the  primary formulation of the Poisson
equation. The reason is that only the stress variable, which can be decoupled from
the primary variable, needs to be involved in the analysis.
This is not,  however,  the case for the Stokes problem under
consideration here because the two variables, velocity and pressure,
are coupled and cannot be separated in  analyses of the convergence and optimality.  To circumvent this difficulty,
B\"{a}nsch, Morin, and Nochetto developed a modified adaptive
procedure in which the Uzawa algorithm on the continuous level is used as
the outer iteration
\cite{BanschMorinNochetto02,Kondratyuk06,KondratyukStevenson07}.

The optimality of the adaptive finite element method for the Poisson equation  is analyzed based on
discrete reliability (see \cite{CasconKreuzerNochettoSiebert07,Stevenson06,Stevenson05} and the references therein).
Basically,  we need one restriction operator
and one prolongation operator in order to analyze the discrete reliability.   For the conforming method,  a  natural  candidate
for the prolongation operator is the usual inclusion operator, and  for the restriction operator  a Scott--Zhang--type
  can be used as it has both the local projection property and the  global and uniform boundedness  property. For  the nonconforming method under consideration here, however,
it is a challenge to come up with a  prolongation operator that has both the local projection property and the  global and uniform boundedness property. For the nonconforming linear element method for the Poisson equation,  such a difficulty can be circumvented
using the discrete Helmholtz decomposition \cite{BeckerMaoShi2010,Rabaus2010}. However,  the Helmholtz decomposition  seems not applicable  for the problem under consideration because the existence of such a decomposition is unclear  for the general case.

The first convergence and optimality analysis of a standard adaptive finite element
method for the Stokes problem was presented  in a technical
report~\cite{HuXu2007} in 2007 by the authors of this paper.  The analysis was based on some special relation
between the nonconforming $P_1$ element and the lowest Raviart--Thomas
element for the Stokes problem and one prolongation operator between
the discrete spaces.
 But we later  found  a gap in our  discrete reliability analysis
 caused by the prolongation operator used therein.
  A  convergence and
optimality analysis was published in \cite{BeckerMao2011} in 2011; however, we also found
 a gap   in their analysis similar to that
 in our earlier report \cite{HuXu2007} (see Appendix A for more details).

The present paper is an improved version of~\cite{HuXu2007} with simplified
and corrected
 proofs.  Its  purpose  is to provide a
rigorous analysis of the convergence and optimality of the adaptive
nonconforming linear element method for the Stokes problem.  The main
idea is to establish the orthogonality or quasi--orthogonality of both
the velocity variable and the pressure variable.  The nonconformity of the discrete velocity
space is the main difficulty in  establishing  the desired  quasi--orthogonality  property  and
the discrete reliability estimate.
To overcome this difficulty we take two steps, (1)  we establish
the quasi--orthogonality for both the velocity and pressure variables by  using  a special conservative
property of the nonconforming linear element,  and  (2) we introduce a new prolongation operator that has
both the projection property and the uniform boundedness property
 for  the discrete reliability analysis.  To
analyze  optimality within the standard nonlinear approximate class
\cite{CasconKreuzerNochettoSiebert07}, we define a new interpolation
operator to bound the consistency error and prove that  the consistency error can be bounded by the approximation
error up to oscillation.  This in fact implies  that the nonlinear approximate class used in \cite{HuXu2007} is {\em equivalent}
 to the standard nonlinear approximate class
\cite{CasconKreuzerNochettoSiebert07}.   Finally, by introducing a new parameter-dependent error
estimator, we prove  convergence and optimality estimates for the
Stokes problem.

The rest of the paper is organized as follows.  In Section 2 we
present the Stokes problem and its nonconforming linear finite element
method, and recall a posteriori error estimate according  to
\cite{ca04,CarstensenFunken01a,CarstensenHu07,DariDuranPadra95}.  We prove the quasi--orthogonality in
Section 3 and then show the
reduction of some total error in Section 4 in terms of  a new
parameter-dependent estimator. We introduce a new prolongation operator to
establish  discrete reliability in Section 5. And, we show
optimality of the adaptive nonconforming linear element method in
Section 6.

\section{The adaptive nonconforming linear element}

Let us first introduce some notations. We use the standard gradient and divergence
operators $\na r:=(\pa r/\pa x\,,\pa r/\pa y)$ for a scalar function $r$, and  $
\div\bpsi:={\pa\psi_1}/{\pa x}+{\pa\psi_2}/{\pa y}$ for a vector
function $\bpsi=(\psi_1,\psi_2)$.  Given a polygonal
domain $\Omega\subset \R^2$ with the boundary $\pa\Omega$, we  use
the standard notation for Sobolev spaces, such as $H^1(\Om)$ and
$L^2(\Om)$.  We define
\begin{equation*}
\begin{split}
H_0^1(\Omega):=\{v\in H^1(\Om), v=0 \text{ on } \pa \Om\}\,, \text{ and }\\[1.0ex]
L_0^2(\Omega):=\{q\in L^2(\Om), \int_{\Om}qdx=0\}.
\end{split}
\end{equation*}
In addition, we denote $(\cdot, \cdot)_{L^2(\Om)}$ as the
usual $L^2$ inner product of functions in the space $L^2(\Om)$, and
$\|\cdot\|_{L^2(\Om)}$ the $L^2$ norm.

Suppose that $\overline{\Omega}$ is covered exactly by a  sequence of
shape--regular triangulations $\cT_k$ ($k\geq 0$) consisting of triangles in
$2D$ (see \cite{CiaBook}), and that this sequence is produced by some
adaptive algorithm where $\cT_k$ is some nested refinement of
$\cT_{k-1}$ by the newest vertex bisection \cite{Stevenson06,Stevenson05}. Let $\cE_k$ be the set of
all edges in $\cT_k$; $\cE_k(\Om)$  the set of interior
edges; $\cE(K)$  the set of edges of any given element $K$ in
$\cT_k$; and $h_K=|K|^{1/2}$ the size of the element $K\in
\cT_k$ where $|K|$ is the area of element $K$.  $\omega_K$ is the
union of elements $K'\in \cT_k$ that share an edge with $K$, and
$\omega_E$ is the union of elements that share a common edge
$E$. Given any edge $E\in\cE_k(\Om)$ with the length $h_E$,  we
assign one fixed unit normal $\nu_E:=(\nu_1,\, \nu_2)$ and tangential
vector $\tau_E:=(-\nu_2,\,\nu_1)$. For $E$ on the boundary,  we choose
$\nu_E:=\nu$,  the unit outward normal to $\Omega$.  Once $\nu_E$ and
$\tau_E$ are fixed on $E$, in relation to $\nu_E$ we define
the elements $K_{-}\in \cT_k$ and $K_{+}\in \cT_k$, with
$E=K_{+}\cap K_{-}$.  Given $E\in\cE_k(\Omega)$ and some
$\R^d$-valued function $v$ defined in $\Omega$, with $d=1,2$, we
denote  $[v]:=(v|_{K_+})|_E-(v|_{K_-})|_E$ as the jump of $v$ across
$E$, where $v|_K$ is the restriction of $v$ on $K$ and $v|_E$ is the
restriction of $v$ on $E$.

\subsection{The Stokes problem and its nonconforming linear element }
The Stokes problem is defined as follows: Given $g\in L^2(\Om)^2$, find $(u, p)\in V\x
Q:=(H_0^1(\Om))^2\x L_0^2(\Om)$ such that
\begin{equation}\label{eq7.1}
\begin{split}
&a(u,v)+b(v,p)+b(u,q)=(g,v)_{L^2(\Om)} \text{ for any }(v,q)\in V\x
Q\,,
\end{split}
\end{equation}
where $u$ and $p$ are the velocity and pressure of the flow, respectively,  and
\begin{equation}
a(u,v):=\mu(\na u,\na v)_{L^2(\Om)} \text{ and } b(v,q):=(\div v,
q)_{L^2(\Om)},
\end{equation}
where $\mu>0$ is the viscosity coefficient of the flow.

Given $\omega\subset \R^2$ and some integer $\ell$,  denote $P_{\ell}(\omega)$ as the space
of polynomials of degree $\leq \ell$ over $\omega$.  We define
\begin{equation*}
\begin{split}
V_k:&=\{v_k\in L^2(\Om)^2, v_k|_{K}\in P_1(K)^2 \text{ for any }K\in
\cT_k, \int_E[v_k]\ds=0\\[1.0ex]
&\qquad\text{ for any }E\in\cE_k(\Om), \text{ and } \int_Ev_k\ds=0
\text{
for any }E\in\cE_k\cap\pa\Om\}\,,\\[1.0ex]
 Q_k:&=\{q_k\in Q, q_k|_K\in P_0(K) \text{ for any } K\in \cT_k\}.
\end{split}
\end{equation*}
Since $V_k$ is not a subspace of $H^1(\Om)^2$, the gradient and
divergence operators are  defined element by element with respect to $\cT_k$, and denoted by
$\na_k$ and $\div_k$.
 Define the
piecewise smooth space
\begin{equation}
H^1(\cT_k):=\{v\in L^2(\Om), v|_K\in H^1(K) \text{ for any
}K\in\cT_k\}\,.
\end{equation}
The discrete bilinear forms  read
 \begin{equation}
 a_k(u,v):=\mu(\na_ku,\na_kv)_{L^2(\Om)} \text{ and } b_k(v, q):=(\div_kv, q)_{L^2(\Om)}
 \end{equation}
 for any  $u, v\in (H^1(\cT_k))^2,  \text{ and } q\in Q$.

The nonconforming $P_1$ element,  proposed in \cite{CroRav73},  for
the Stokes problem is as follows: Given $g\in
L^2(\Om)^2$,  find $(u_k,p_k)\in V_k\x Q_k$ such that
\begin{equation}\label{eq7.3}
\begin{split}
\begin{split}
&a_k(u_k,v)+b_k(v,p_k)+b_k(u_k,q)=(g,v)_{L^2(\Om)} \text{ for any
}(v,q) \in V_k\x Q_k\,.
\end{split}
\end{split}
\end{equation}
Let $\id\in \R^{2\x2}$ be the identity matrix. Define
$$
 \sigma_{k}:=\mu\na_ku_k+p_k\id.
$$
Then, we have
\begin{equation}
(\sigma_k, \na_kv_k)_{L^2(\Om)}=(g, v_k)_{L^2(\Om)} \text{ for any }v_k\in V_k.
\end{equation}
\subsection{The a posteriori error estimate}  To recall the a posteriori error estimator of  the nonconforming $P_1$  element,  we define
the residual  $\res_{k-1}(\cdot)$  by
\begin{equation}\label{Res4}
\res_{k-1}(v):=(g,v)_{L^2(\Om)}-a_k(u_{k-1},v)-b_k(v,p_{k-1}) \text{ for any
}v\in H^1(\cT_k)^2\,,
\end{equation}
with the solution $(u_{k-1}, p_{k-1})$ of \eqref{eq7.3} on the mesh $\cT_{k-1}$,
which is a coarser and nested mesh of $\cT_k$.  It follows from the definition of $(u_{k-1}, p_{k-1})$  that
$$
\res_{k-1}(v_{k-1})=0 \text{ for any }v_{k-1}\in V_{k-1}.
$$

Given $K\in\cT_k$, we define the element estimator
\begin{equation}
\eta_{K}(u_k, p_k):=h_K\|g\|_{L^2(K)}+(\sum\limits_{E\subset\pa
K}h_K\|[\na_ku_k  \tau_E]\|_{L^2(E)}^2)^{1/2}.
\end{equation}
Given $S_k\subset\cT_k$, we define the  estimator over it by
\begin{equation}\label{eq7.12}
\eta^2(u_k, p_k, S_k):= \sum\limits_{K\in S_k}\eta_{K}^2(u_k, p_k).
\end{equation}
Given any $K\in \cT_k$, denote $g_K$ as the $L^2$ projection of $g$
onto $P_{0}(K)$. We define the oscillation
\begin{equation}\label{eq5.1b}
\osc^2(g,\cT_k):=\sum\limits_{K\in \cT_k}h_K^2\|g-g_K\|_{L^2(K)}^2.
\end{equation}

The  reliability and efficiency of the estimator $\eta(u_k,
p_k, \cT_k)$ can be found  in
\cite{ca04,CarstensenFunken01a,CarstensenHu07,DariDuranPadra95}, as stated in the following lemma.

\begin{lemma}\label{Theorem7.2} Let $(u,p)$ and $(u_k, p_k)$ be the solutions
of the Stokes problem
 \eqref{eq7.1} and the discrete problem  \eqref{eq7.3}, respectively. Then,
\begin{equation}\label{eq7.11}
\|\na_k(u-u_k)\|_{L^2(\Om)}^2+\|p-p_k\|_{L^2(\Om)}^2\lesssim
\eta^2(u_k, p_k, \cT_k),
\end{equation}
\begin{equation}\label{eq7.11b}
\eta^2(u_k, p_k, \cT_k)\lesssim
\|\na_k(u-u_k)\|_{L^2(\Om)}^2+\|p-p_k\|_{L^2(\Om)}^2+\osc^2(g,
\cT_k).
\end{equation}
\end{lemma}
\begin{remark}
For the Stokes problem, the estimator usually involves the pressure approximation. For
the nonconforming $P_1$ element, as shown in the above lemma,  we can decouple
the pressure from the velocity \cite{DariDuranPadra95}.
\end{remark}
Here and throughout the paper, we  use the notations $\lesssim$
and $\approxeq$. When we write
$$
A_1 \lesssim B_1, \text{ and } A_2\approxeq B_2,
$$
 possible constants $C_1$, $c_2$ and $C_2$ exist such that
$$
A_1 \leq C_1 B_1, \text{ and } c_2B_2\leq A_2\leq C_2 B_2.
$$
\subsection{The adaptive nonconforming finite element method}
The adaptive  algorithm  is defined as follows: Let
$\mathcal T_0$ be an  initial shape--regular triangulation, a right--side $g\in L^2(\Omega)^2$, a tolerance $\epsilon$, and a parameter  $0<\theta<1$.

\begin{algorithm}\label{Algorithm}
\noindent $[\mathcal{T}_N\,, u_N\,, p_N]$={\bf \small
ANFEM}$(\mathcal{T}_0, g, \epsilon, \theta)$

$\eta = \epsilon\,, k=0$

\smallskip

{\bf \small WHILE} $\eta \geq \epsilon$, {\bf \small DO}
 \begin{enumerate}

\item Solve \eqref{eq7.3} on $\mathcal{T}_{k}$ to get the solution
$(u_{k}, p_{k})$.

 \item Compute the error estimator $\eta=\eta(u_{k}, p_{k}, \cT_{k})$.

\item Mark the minimal element set $\mathcal{M}_{k}$ such that
  \begin{equation}\label{bulk}
 \eta ^2(u_{k}, p_{k}, \mathcal{M}_{k})\geq \theta \, \eta ^2(u_{k}, p_{k}, \cT_{k}).
  \end{equation}

  \item Refine each triangle $K \in \mathcal{M}_{k}$
  by the newest vertex bisection to get $\cT_{k+1}$ and set $k=:k+1$.
\end{enumerate}

{\bf \small END WHILE}

\smallskip

$\cT_N=\cT_k$.

\noindent {\bf \small END ANFEM}
\end{algorithm}

\section{Quasi--orthogonality}
 The quasi--orthogonality property is the main
ingredient for  the convergence  analysis of the adaptive
nonconforming method under consideration.  In this section we  establish  such a property by
exploring the conservative property of the nonconforming linear
element  and by confirming  that the stress is piecewise constant. To this
end, we  define a canonical interpolation operator $\Pi_k$ for the
nonconforming space $V_k$
 and a restriction operator $I_{k-1}$ from $V_k$ to the coarser space $V_{k-1}$.
 Given $v\in V$, we define  the interpolation $\Pi_kv\in V_k$ by
 \begin{equation}\label{interpolation}
 \int_E \Pi_kvds:=\int_E vds\text{  for any }E\in\cE_k\,.
 \end{equation}
 In this paper, the above property is referred to  as the conservative property. This property is crucial for the analysis herein.
  A similar conservative property was first explored in \cite{HSX10} to analyze  the quasi--orthogonality property  of the Morley
  element.

The  interpolation admits the following estimate:
\begin{equation}\label{eq3.4a}
\|v-\Pi_kv\|_{L^2(K)}\lesssim h_K\|\nabla v\|_{L^2(K)}\text{ for any
}K\in\cT_k\text{ and }v\in V\,.
\end{equation}
 Given $v_k\in V_k$, we define the restriction interpolation $I_{k-1}v_k\in V_{k-1}$  by
 \begin{equation}\label{Restriction}
\int_E I_{k-1}v_k ds:=\sum\limits_{l=1}^{\ell}\int_{E_l} v_k ds\,,  E\in
\cE_{k-1} \text{ with }E=E_1\cup E_2\cdots\cup E_{\ell}\text{
 and }E_i\in\cE_k\,.
\end{equation}
The properties of the restriction operator $I_{k-1}$ are summarized in
the following lemma.
\begin{lemma} Let the restriction operator $I_{k-1}$  be defined
in \eqref{Restriction}. Then,
\begin{equation}\label{eq3.4b}
I_{k-1}v_k=v_k  \text{ for any } K\in \cT_k\cap\cT_{k-1},  v_k\in V_k\,,
\end{equation}
\begin{equation}\label{eq3.5b}
\|I_{k-1}v_k-v_k\|_{L^2(K)}\lesssim h_K\|\nabla_kv_k\|_{L^2(K)} \text{ for any
}K\in\cT_{k-1}\backslash\cT_k,  v_k\in V_k\,.
\end{equation}
\end{lemma}
\begin{proof} The property  \eqref{eq3.4b} directly follows  from  the definition of the
restriction interpolation.  Only the estimate
\eqref{eq3.5b} needs to be proved.
  In fact, both sides of \eqref{eq3.5b}
are semi-norms of the restriction $(V_k)_K$ of $V_k$ on $K$. If the
right--hand side vanishes for some $v\in (V_k)_K$, then $v_k$ is a
piecewise constant vector over $K$ with respect to $\cT_k$. Given the average continuity of $v_k$ across the internal edges of $\cT_k$,
 it follows that $v_k$ is a constant vector  on $K$. Therefore, the left--hand side also
vanishes for the same $v_k$. The desired result then follows  a
scaling argument.
\end{proof}
\begin{remark} An alternative proof for the inequality \eqref{eq3.5b} follows  the discrete Poincare inequality
 established in \cite{Brenner2003} for the scalar function, which is further investigated in \cite{Rabaus2010}.
 Notice that the positive constant of \eqref{eq3.5b} is independent of the ratio
 \begin{equation}
\gamma:=\max\limits_{K\in\cT_{k-1}\backslash\cT_{k}}\max\limits_{\cT_{k}\ni T\subset K}\frac{h_K}{h_T},
\end{equation}
 see \cite[Lemma 4.1]{Rabaus2010}
 for more details.
\end{remark}

\begin{lemma}  Let $(u_{k-1}, p_{k-1})$ be the solution of the discrete
problem \eqref{eq7.3} on the mesh $\cT_{k-1}$.  It, therefore, holds that
\begin{equation}\label{estimateresidual}
|\res_{k-1}(v_k)|\lesssim
\big(\sum\limits_{K\in\cT_{k-1}\backslash\cT_k}h_K^2\|g\|_{L^2(K)}^2\big)^{1/2}\|\na_kv_k\|_{L^2(\Om)}\text{
for any }v_k\in V_k.
\end{equation}
\end{lemma}
\begin{proof} For the reader's convenience,  we recall the definition of the
residual as follows:
\begin{equation}\label{term-1}
\res_{k-1}(v_k)=(g,v_k)_{L^2(\Om)}-(\sigma_{k-1},\na_kv_k)_{L^2(\Om)}.
\end{equation}
To analyze the right-hand side of the above equation, we set
$v_{k-1}=I_{k-1}v_k$.  As $\sigma_{k-1}$ is a piecewise constant tensor
with respect to the mesh $\cT_{k-1}$,  the definition of the
interpolation operator $I_{k-1}$ in \eqref{Restriction} leads to
\begin{equation}\label{term1}
\int_E (v_k-v_{k-1})\cdot \sigma_{k-1}\nu_Eds=0\text{ for any }E\in \cE_{k-1}\,.
\end{equation}
For any $E\in \cE_k$ that lies in the interior of some $K\in\cT_{k-1}$,
 the integral average of $v_k$ over
$E$ is continuous and  $\sigma_{k-1}$ is a constant on $K$.   Then,
\begin{equation}\label{term2}
\int_E [v_k-v_{k-1}]\cdot \sigma_{k-1}\nu_Eds=0.
\end{equation}
By integrating  parts on the fine mesh $\cT_k$ and using \eqref{term1} and \eqref{term2},  we get
\begin{equation}
(\na_k(v_k-v_{k-1}),\sigma_{k-1})_{L^2(\Om)}=0.
\end{equation}
Inserting this identity into \eqref{term-1}  and  adopting the discrete
problem \eqref{eq7.3},  we employ   properties
\eqref{eq3.4b} and \eqref{eq3.5b} of the interpolation operator
$I_{k-1}$ to derive
\begin{equation}
\begin{split}
|\res_{k-1}(v_k)|&=|(g,v_k-v_{k-1})_{L^2(\Om)}|
\leq\sum\limits_{K\in\cT_{k-1}\backslash\cT_k}\|g\|_{L^2(K)}\|v_k-v_{k-1}\|_{L^2(K)}\\[0.5ex]
&\lesssim
\sum\limits_{K\in\cT_{k-1}\backslash\cT_k}h_K\|g\|_{L^2(K)}\|\na_kv_k\|_{L^2(K)}\,,
\end{split}
\end{equation}
which completes the proof.
\end{proof}

\begin{lemma}(Quasi-orthogonality of the velocity) \label{Lemma7.5}
Let $(u_k,p_k)$ and $(u_{k-1}, p_{k-1})$ be the discrete
 solutions of \eqref{eq7.3} on $\cT_k$ and $\cT_{k-1}$, respectively.
 Then,
\begin{equation*}
\begin{split}
|a_k(u-u_k,u_k-u_{k-1})|\lesssim \|\na_k(u-u_k)\|_{L^2(\Om)}
\bigg(\sum\limits_{K\in\cT_{k-1}\backslash\cT_k}h_K^2\|g\|_{L^2(K)}^2\bigg)^{1/2}\,.
 \end{split}
\end{equation*}
\end{lemma}
\begin{proof} The Stokes problem \eqref{eq7.1} and  the discrete problem \eqref{eq7.3} give
\begin{equation}\label{eq7.26}
\begin{split}
a_k(u-u_k,u_k-u_{k-1}) =(\na_k(u-u_k), \sigma_k-\sigma_{k-1})_{L^2(\Om)}.
\end{split}
\end{equation}
Given that $(\div_k(u-u_k), p_k-p_{k-1})_{L^2(\Om)}=0$, let $v_k=\Pi_k(u-u_k)$.
 And,  $\sigma_k-\sigma_{k-1}$ is a piecewise
constant tensor with respect to the fine mesh  $\cT_k$; therefore,  by the definition
of the interpolation operator $\Pi_k$ in \eqref{interpolation}, we
integrate by parts on $\cT_k$ to obtain
\begin{equation}
(\na_k((u-u_k)-v_k), \sigma_k-\sigma_{k-1})_{L^2(\Om)}=0.
\end{equation}
From the discrete problem \eqref{eq7.3},  we have
\begin{equation}\label{term0}
a_k(u-u_k,u_k-u_{k-1})=(g,v_k)_{L^2(\Om)}-(\na_kv_k,\sigma_{k-1})_{L^2(\Om)}=\res_{k-1}(v_k).
\end{equation}
The term on the right-hand side of the equation \eqref{term0} can be
estimated by the inequality \eqref{estimateresidual} as follows:
\begin{equation*}
\begin{split}
|\res_{k-1}(v_k)| &\lesssim
\sum\limits_{K\in\cT_{k-1}\backslash\cT_k}h_K\|g\|_{L^2(K)}\|\na_kv_k\|_{L^2(K)}\\[0.5ex]
&\lesssim
\sum\limits_{K\in\cT_{k-1}\backslash\cT_k}h_K\|g\|_{L^2(K)}\|\na_k(u-u_k)\|_{L^2(K)}\,,
\end{split}
\end{equation*}
which completes the proof.
\end{proof}

\begin{lemma}(Quasi--orthogonality of the pressure)\label{Lemma7.6}
Let $(u_k,p_k)$ and $(u_{k-1}, p_{k-1})$ be the discrete
 solutions of \eqref{eq7.3} on $\cT_k$ and $\cT_{k-1}$,  respectively.
 Then,
 \begin{equation}\label{eq7.36}
 \begin{split}
&|(p-p_k, p_k-p_{k-1})_{L^2(\Om)}|\\[0.3ex]
 & \lesssim
 \bigg(\big(\sum\limits_{K\in\cT_{k-1}\backslash\cT_k}
h_K^2\|g\|_{L^2(K)}^2\big)^{1/2}
 +\|\na_k(u_k-u_{k-1})\|_{L^2(\Om)}\bigg)\|p-p_k\|_{L^2(\Om)}\,.
 \end{split}
 \end{equation}
\end{lemma}
\begin{remark} The quasi--orthogonality of  the pressure herein is different
from those for the nonstandard method of the Poisson equation
\cite{CarstensenHoppe05a,CarstensenHoppe05b,ChenHolstXu07} by the fact that both
$\|\na_k(u_k-u_{k-1})\|_{L^2(\Om)}$ and $\|p-p_k\|_{L^2(\Om)}$ appear on the right--hand side of \eqref{eq7.36}.
\end{remark}

\begin{proof}  Let $\Pi_{0, k}$ be the $L^2$ projection operator from $L^2_0(\Omega)$ onto $Q_k$.
 It follows from the discrete inf-sup condition that there exists $v_k\in V_k$ with
\begin{equation}\label{eq7.37}
\div_k v_k=\Pi_{0, k}p-p_k, \text{ and
}\|\na_kv_k\|_{L^2(\Om)}\lesssim \|\Pi_{0, k}p-p_k\|_{L^2(\Om)}.
\end{equation}
Since $p_k-p_{k-1}\in Q_k$, it follows from  the continuous problem
\eqref{eq7.1},
 the discrete problem \eqref{eq7.3}, and the definition of the residual \eqref{Res4} that
\begin{equation*}
\begin{split}
(p-p_k, p_k-p_{k-1})_{L^2(\Om)}
 =(\div_k v_k, p_k-p_{k-1})_{L^2(\Om)}
 =\res_{k-1}(v_k)+a_k(u_{k-1}-u_k,v_k).
 \end{split}
\end{equation*}
We use the estimates in  \eqref{estimateresidual} and \eqref{eq7.37}
to get
\begin{equation*}
\begin{split}
&|(p-p_k, p_k-p_{k-1})_{L^2(\Om)}|\\[0.5ex]
& \lesssim \bigg(\big(\sum\limits_{K\in\cT_{k-1}\backslash\cT_k}
h_K^2\|g\|_{L^2(K)}^2\big)^{1/2}+\|\na_k(u_k-u_{k-1})\|_{L^2(\Omega)}\bigg)\|p-p_k\|_{L^2(\Omega)},
\end{split}
\end{equation*}
which completes the proof.
\end{proof}

\section{The convergence of the ANFEM}
To prove the  convergence of the adaptive algorithm, we first
prove the reduction of the error between  the two nested  meshes,  $\cT_k$
and $\cT_{k-1}$,  where $\cT_k$ is the refinement of the coarser mesh
$\cT_{k-1}$ with \eqref{bulk}  by the newest vertex bisection.  In order to control the volume part
$\sum\limits_{K\in\cT_{k-1}\backslash\cT_k}h_K^2\|g\|_{L^2(K)}^2$  appearing in Lemmas \ref{Lemma7.5} and \ref{Lemma7.6}, we introduce the following modified estimator:
\begin{equation}\label{eq7.53b}
\begin{split}
&\tilde{\eta}^2(u_{k-1}, p_{k-1},
\cT_{k-1}):=\sum\limits_{K\in\cT_{k-1}}\big(\beta_1
h_K^2\|g\|_{L^2(K)}^2+\eta_K^2(u_{k-1}, p_{k-1})\big)\,
\end{split}
\end{equation}
with the positive constant $\beta_1>0$ to be determined later. Note
that this modified estimator is introduced only for the  convergence  analysis and that
 the final convergence and optimal complexity will be proved for Algorithm \ref{Algorithm}.

Note that the volume  residual
$\sum\limits_{K\in\cT_{k-1}} h_K^2\|g\|_{L^2(K)}^2$ does not contain the
unknowns.  Hence, we add it to settle down the
lacking of  the Galerkin--orthogonality or quasi--orthogonality.  We stress that
 the Galerkin--orthogonality or quasi--orthogonality is an essential ingredient for the  convergence analysis
 of the adaptive conforming, nonconforming, and mixed methods for the
Poisson-like problems \cite{CarstensenHoppe05a, CarstensenHoppe05b,CasconKreuzerNochettoSiebert07,ChenHolstXu07,Dolfler96,MekchayNochetto05,
MorinNochettoSiebert00}.  This is
another reason that we need a modified estimator as in
\eqref{eq7.53b}.

We list  three standard  components  for the convergence
analysis of the adaptive method, which can be proved by following
the arguments, for  instance, in
\cite{CarstensenHoppe05b,CasconKreuzerNochettoSiebert07,Dolfler96}.

\begin{lemma}\label{Lemma5.3b} Let $\cT_k$ be some refinement of $\cT_{k-1}$ from Algorithm \ref{Algorithm},
 then  $\rho>0$ and a
  positive constant $\beta\in(1-\rho\theta,1)$ exist,
such that
\begin{equation}\label{eq5.7}
{\eta}^2(u_{k-1},p_{k-1},\cT_k)
 \leq  \beta {\eta}^2(u_{k-1},p_{k-1},\cT_{k-1})+(1-\rho\theta-\beta){\eta}^2(u_{k-1}, p_{k-1}, \cT_{k-1})\,.
\end{equation}\end{lemma}
\begin{proof} The result can be proved by following the idea in \cite{CarstensenHoppe05b,CasconKreuzerNochettoSiebert07,Dolfler96}.
 The details are only given  for the readers' convenience.  In fact, we have
\begin{equation}
{\eta}^2(u_{k-1}, p_{k-1}, \cT_k)={\eta}^2(u_{k-1}, p_{k-1}, \cT_{k-1}\cap\cT_k)+{\eta}^2(u_{k-1}, p_{k-1}, \cT_k\backslash\cT_{k-1}).
\end{equation}
For any $K\in \cT_{k-1}\backslash\cT_k$, we only need to consider the
case where  $K$ is subdivided into $K_{1}\,,K_2\in\cT_k$ with
$|K_1|=|K_2|=\frac{1}{2}|K|$.  As $[\nabla_{k-1}u_{k-1} \tau_E]=0$ over the interior edge $E=K_1\cap K_2\in\cE_k$,   we  have
\begin{equation}\label{eq4.5}
\begin{split}
 &\sum\limits_{i=1}^2{\eta}_{K_i}^2(u_{k-1}, p_{k-1})\\
& :=\sum\limits_{i=1}^2\bigg(h_{K_i}\|g\|_{L^2(K_i)}+\bigg(\sum\limits_{\cE_k\ni
E\subset\pa
K_i}h_{K_i}\|[\nabla_{k-1}u_{k-1} \tau_E]\|_{L^2(E)}^2\bigg)^{1/2}\bigg)^2\\[0.5ex]
& \leq\frac{1}{2^{1/2}} {\eta}^2_K(u_{k-1}, p_{k-1})\\
&:=\frac{1}{2^{1/2}}
 \bigg( h_{K}\|g\|_{L^2(K)}+\bigg(\sum\limits_{\cE_{k-1}\ni
E\subset\pa
K}h_{K}\|[\nabla_{k-1}u_{k-1}\tau_E]\|_{L^2(E)}^2\bigg)^{1/2}\bigg)^2\,.
 \end{split}
\end{equation}
Consequently,
\begin{equation}
\begin{split}
\sum\limits_{K\in\cT_{k-1}\backslash\cT_k}
\sum\limits_{i=1}^2{\eta}_{K_i}^2(u_{k-1}, p_{k-1}) \leq
\frac{1}{2^{1/2}}{\eta}^2(u_{k-1}, p_{k-1}, \cT_{k-1}\backslash\cT_k)\,.
\end{split}
\end{equation}
Let $\rho=1-\frac{1}{2^{1/2}}$, therefore,  we obtain
\begin{equation}
\begin{split}
{\eta}^2(u_{k-1}, p_{k-1}, \cT_k) &\leq {\eta}^2(u_{k-1}, p_{k-1}, \cT_{k-1})-\rho{\eta}^2(u_{k-1}, p_{k-1},
\cT_{k-1}\backslash\cT_k)\,.
\end{split}
\end{equation}
 Choosing the positive parameter
$\beta$ with $1-\rho\theta<\beta<1$,  we combine the above inequality and the bulk criterion \eqref{bulk} to
 achieve  the desired result.
\end{proof}

\begin{lemma}\label{rightcon}
Let $\cT_k$ be some refinement of $\cT_{k-1}$  produced in Algorithm \ref{Algorithm},
 then there exists $\rho>0$  such that
 \begin{equation}\label{righthandside}
 \sum\limits_{K\in\cT_k}h_K^{2}\|g\|_{L^2(K)}^2\leq
 \sum\limits_{K\in\cT_{k-1}}h_K^{2}\|g\|_{L^2(K)}^2-\rho\sum\limits_{K\in\cT_{k-1}\backslash\cT_k}h_K^{2}\|g\|_{L^2(K)}^2\,.
 \end{equation}
\end{lemma}
\begin{proof}  This can be proved by a similar argument proposed in the previous lemma.
\end{proof}

\begin{lemma}\label{continuity}(Continuity of the estimator)  Let $u_k$ and $u_{k-1}$ be the solutions to the discrete problem \eqref{eq7.3}
on the meshes $\cT_k$ and $\cT_{k-1}$  obtained from Algorithm \ref{Algorithm}.  Given any positive
constant $\epsilon$, there exists a positive constant
$\beta_2(\epsilon)$ dependent on $\epsilon$ such that
\begin{equation}\label{eq5.4}
{\eta}^2(u_k, p_k,  \cT_k)\leq (1+\epsilon) {\eta}^2(u_{k-1}, p_{k-1},
\cT_k)
+\frac{1}{\beta_2(\epsilon)}\|\na_k(u_k-u_{k-1})\|_{L^2(\Om)}^2\,.
\end{equation}
\end{lemma}
\begin{proof}
Given any $K\in\cT_k$,  it follows from the definitions of
${\eta}_K(u_k, p_k)$ and ${\eta}_K(u_{k-1}, p_{k-1})$ in \eqref{eq4.5} that
\begin{equation*}
\begin{split}
&\big |{\eta}_K(u_k, p_k)-{\eta}_K(u_{k-1}, p_{k-1})\big |\\
&=\bigg
|\bigg(\sum\limits_{\cE_k\ni E\subset\pa K}h_K\|[\nabla_k
 u_k\tau_E]\|_{L^2(E)}^2\bigg)^{1/2}
 -\bigg(\sum\limits_{\cE_k\ni
E \subset\pa K}h_K\|[\nabla_{k-1}
 u_{k-1}\tau_E]\|_{L^2(E)}^2\bigg)^{1/2}\bigg |\\
 &\leq \bigg(\sum\limits_{\cE_k\ni
E \subset\pa K}h_K\|[\nabla_k
 (u_k-u_{k-1})\tau_E]\|_{L^2(E)}^2 \bigg)^{1/2}.
 \end{split}
\end{equation*}
Given $E\in\cE_k$, let $K_1, K_2\in \cT_k$ be the two elements that take $E$ as one edge. Then, we use the trace theorem and the fact
 that $\nabla_k(u_k-u_{k-1})$ is a piecewise constant tensor to get
\begin{equation}
\begin{split}
&\|[\nabla_k
 (u_k-u_{k-1})\tau_E]\|_{L^2(E)}\\
 &\leq \|\nabla_k
 (u_k-u_{k-1})\tau_E|_{K_1}\|_{L^2(E)}+\|\nabla_k
 (u_k-u_{k-1})\tau_E|_{K_2}\|_{L^2(E)}\\[0.5ex]
 &\lesssim h_K^{-1/2}\|\nabla_k
 (u_k-u_{k-1})\|_{L^2(\omega_E)}\,,
 \end{split}
\end{equation}
which gives
\begin{equation}
{\eta}_K(u_k, p_k)\leq {\eta}_K(u_{k-1}, p_{k-1})+ C_{Con} \|\nabla_k
 (u_k-u_{k-1})\|_{L^2(\omega_K)},
\end{equation}
for some positive constant $C_{Con}$.  Given any positive constant $\epsilon$, we apply the Young inequality to get
\begin{equation}
{\eta}^2_K(u_k, p_k)\leq (1+\epsilon) {\eta}_K^2(u_{k-1}, p_{k-1})+ \frac{C_{Con}^2(1+\epsilon)}{\epsilon} \|\nabla_k
 (u_k-u_{k-1})\|_{L^2(\omega_K)}^2.
\end{equation}
A summation  over all elements in $\cT_k$ completes the proof with  $\beta_2(\epsilon)=\frac{M\epsilon}{C_{Con}^2(1+\epsilon)}$,
 where the positive constant $M$ depends on the finite overlapping of the patches $\omega_K$.
\end{proof}

In the following theorem,  we prove the convergence of the adaptive
nonconforming finite element method for the  Stokes problem.  The main ingredients are the quasi--orthogonality
 of both the velocity and the pressure in Lemmas \ref{Lemma7.5} and \ref{Lemma7.6},  and  the relations of the estimators
  between two the meshes  $\cT_k$ and $\cT_{k-1}$ presented in Lemmas \ref{Lemma5.3b}--\ref{continuity}.
\begin{theorem}\label{Theorem7.8} Let $(u,p)$ and $(u_k,p_k)$
 be  the solutions of \eqref{eq7.1} and \eqref{eq7.3}.
 Then  $\gamma_1, \gamma_2, \beta_1>0$ and  $0<\al<1$ exist,  such that
\begin{equation}\label{eq7.44}
\begin{split}
&\|\na_k(u-u_k)\|_{L^2(\Om)}^2+\gamma_1\|p-p_k\|_{L^2(\Om)}^2+\gamma_2\tilde{\eta}^2(u_k,
p_k,
\cT_k)\\[0.5ex]
& \leq \al\big(
\|\na_{k-1}(u-u_{k-1})\|_{L^2(\Om)}^2+\gamma_1\|p-p_{k-1}\|_{L^2(\Om)}^2+\gamma_2\tilde{\eta}^2(u_{k-1},
p_{k-1}, \cT_{k-1})\big).
\end{split}
\end{equation}
\end{theorem}
\begin{proof}
First, we adopt the quasi--orthogonality of  both the velocity and the pressure.  Denote the multiplication constant in Lemma
\ref{Lemma7.5} by $C_{QOV}$.  As
\begin{equation}
\begin{split}
\|\na_k(u-u_k)\|_{L^2(\Om)}^2&=\|\na_k(u-u_{k-1}\|_{L^2(\Om)}^2-\|\na_k(u_k-u_{k-1})\|_{L^2(\Om)}^2\\
&\quad -2(\na_k(u-u_k), \na_k(u_k-u_{k-1}))_{L^2(\Om)},
\end{split}
\end{equation}
  it follows from the  quasi--orthogonality of the velocity  in Lemma
\ref{Lemma7.5}  and  the Young inequality that
\begin{equation}\label{first}
\begin{split}
&(1-\delta_1)\|\na_k(u-u_k)\|_{L^2(\Om)}^2\\
&\leq
\|\na_{k-1}(u-u_{k-1})\|_{L^2(\Om)}^2
-\|\na_k(u_k-u_{k-1})\|_{L^2(\Om)}^2\\
&\quad +C_{1}(\delta_1)\sum\limits_{K\in\cT_{k-1}\backslash\cT_k}h_K^2\|g\|_{L^2(K)}^2,
\end{split}
\end{equation}
where $C_1(\delta_1)=\frac{C_{QOV}^2}{\delta_1}$ for any positive constant $0<\delta_1<1$. Denote the multiplication constant in Lemma \ref{Lemma7.6}
 by $C_{QOP}$.  From the quasi--orthogonality of the pressure proved in Lemma
\ref{Lemma7.6} and the Young inequality,  we have
\begin{equation}\label{second}
\begin{split}
(1-\delta_2-\delta_3)\|p-p_k\|_{L^2(\Om)}^2& \leq
\|p-p_{k-1}\|_{L^2(\Om)}^2-\|p_k-p_{k-1}\|_{L^2(\Om)}^2 \\[1.0ex]
&\quad
+\frac{1}{\beta_3(\delta_3)}\|\na_k(u_k-u_{k-1})\|_{L^2(\Om)}^2\\[1.ex]
&\quad +C_{2}(\delta_2)
\sum\limits_{K\in\cT_{k-1}\backslash\cT_k}h_K^2\|g\|_{L^2(K)}^2;
\end{split}
\end{equation}
here $\beta_3(\delta_3)=\frac{\delta_3}{C_{QOP}^2}$ and
 $C_2(\delta_2)=\frac{C_{QOP}^2}{\delta_2}$ for any constants $0<\delta_2, \delta_3<1$.
 Then we multiply the inequality \eqref{first} by $\gamma_1>0$ and the inequality \eqref{second} by
$\gamma_2>0$ to obtain
\begin{equation}\label{third}
\begin{split}
&\gamma_1(1-\delta_1)\|\na_{k-1}(u-u_k)\|_{L^2(\Om)}^2+\gamma_2(1-\delta_2-\delta_2)\|p-p_k\|_{L^2(\Om)}^2\\[0.5ex]
& \leq
\gamma_1\|\na_{k-1}(u-u_{k-1})\|_{L^2(\Om)}^2+\gamma_2\|p-p_{k-1}\|_{L^2(\Om)}^2
-(\gamma_1-\frac{\gamma_2}{\beta_3(\delta_3)})\|\na_k(
u_k-u_{k-1})\|_{L^2(\Om)}^2\\[0.5ex]
&\quad-\gamma_2\|p_k-p_{k-1}\|_{L^2(\Om)}^2
+\big(\gamma_1C_{1}(\delta_1)+\gamma_2C_{2}(\delta_2)\big)
\sum\limits_{K\in\cT_{k-1}\backslash\cT_k}h_K^2\|g\|_{L^2(K)}^2.
\end{split}
\end{equation}
For the presentation, we introduce some short--hand  notations
for any positive constants $\gamma_3, \gamma_4>0$:
\begin{equation*}
\begin{split}
\mathfrak{G}_k(u_k,p_k):&=\gamma_1(1-\delta_1)\|\na_{k-1}(u-u_k)\|_{L^2(\Om)}^2
+\gamma_2(1-\delta_2-\delta_3)\|p-p_k\|_{L^2(\Om)}^2\\[0.5ex]
&\quad+\gamma_3\eta^2(u_k, p_k, \cT_k)
+\gamma_4\sum\limits_{K\in\cT_k}h_K^2\|g\|_{L^2(K)}^2\,,
\end{split}
\end{equation*}
\begin{equation}
\begin{split}
\overline{\mathfrak{G}}_{k-1}(u_{k-1}, p_{k-1}):&=\gamma_1\|\na_{k-1}(u-u_{k-1})\|_{L^2(\Om)}^2+\gamma_2\|p-p_{k-1}\|_{L^2(\Om)}^2\\[0.5ex]
&\quad+\gamma_3\beta\eta^2(u_{k-1}, p_{k-1}, \cT_{k-1})
+\gamma_4\sum\limits_{K\in\cT_{k-1}}h_K^2\|g\|_{L^2(K)}^2\,.
\end{split}
\end{equation}

Second,  we  use the continuity of the estimators from Lemmas \ref{Lemma5.3b}--\ref{continuity}  to cancel
 both the term  $\|\na_k(u_k-u_{k-1})\|_{L^2(\Om)}$ and the volume estimator.  In fact,
  from \eqref{eq5.7} and \eqref{eq5.4}, we have
  \begin{equation}
\begin{split}
  \eta^2(u_k, p_k, \cT_k)\leq \beta\eta^2(u_{k-1}, p_{k-1}, \cT_{k-1})+\frac{1}{\beta_2(\epsilon)}\|\na_k(u_k-u_{k-1})\|_{L^2(\Om)}^2\\
  +\big((1-\rho\theta-\beta)(1+\epsilon)+\epsilon\beta\big){\eta}^2(u_{k-1},p_{k-1},\cT_{k-1}).
\end{split}
  \end{equation}
 Then we combine the above inequality with  the inequalities \eqref{third}  and \eqref{righthandside}
  to obtain
\begin{equation*}
\begin{split}
\mathfrak{G}_k(u_k,p_k) &\leq \overline{\mathfrak{G}}_{k-1}(u_{k-1},p_{k-1})
-\big(\gamma_1-\frac{\gamma_2}{\beta_3(\delta_3)}-\frac{\gamma_3}{\beta_2(\epsilon)}\big)\|\na_k(
u_k-u_{k-1})\|_{L^2(\Om)}^2\\[0.5ex]
&\quad-\gamma_2\|p_k-p_{k-1}\|_{L^2(\Om)}^2+\gamma_3\big((1-\rho\theta-\beta)(1+\epsilon)+\epsilon\beta\big){\eta}^2(u_{k-1},
p_{k-1},
\cT_{k-1})\\[0.5ex]
&\quad+\big(\gamma_1C_{1}(\delta_1)+\gamma_2C_{2}(\delta_2)-\gamma_4\rho\big)
\sum\limits_{K\in\cT_{k-1}\backslash\cT_k}h_K^2\|g\|_{L^2(K)}^2.
\end{split}
\end{equation*}
It remains to prove that the  positive constants $\delta_i, i=1, 2, 3$, $\gamma_i, i=1, 2, 3, 4$, $\epsilon$,   $\beta$, and $\beta_1$ exist such that the contraction \eqref{eq7.44} holds for some constant $0<\al<1$. Further it is possible that the constant dependent on
the choices of the aforementioned parameters but independent of the meshsize $h$ and the level $k$.  This will be achieved in the following three steps.

 {\bf Step 1~}  For the second,  fourth, and  fifth terms on the right-hand side of the above inequality to vanish, we set
 \begin{equation}
 \begin{split}
\gamma_2&=(\gamma_1-\frac{\gamma_3}{\beta_2(\epsilon)})\beta_3(\delta_3) \text{ with } \gamma_1>\frac{\gamma_3}{\beta_2(\epsilon)},\\[0.5ex]
\gamma_4&=(\gamma_1C_{1}(\delta_1)+\gamma_2C_{2}(\delta_2))/\rho,\\[0.5ex]
\beta&=(1-\rho\theta)(1+\epsilon).
\end{split}
\end{equation}
Note that $\gamma_2$, $\gamma_4$, and $\beta$ will be  determined after $\delta_i, i=1, 2, 3$,  $\gamma_1$,  $\gamma_3$, and $\epsilon$ have been specified.
In the following, we assume that $\epsilon$ is fixed in such a way  that $0<\beta<1$.
Also,  we let $\gamma_1$ and $\gamma_3$ be fixed such that $\gamma_1>\frac{\gamma_3}{\beta_2(\epsilon)}$ and $\gamma_2>0$.
Hence, we have
\begin{equation*}
\mathfrak{G}_k(u_k,p_k) \leq \overline{\mathfrak{G}}_{k-1}(u_{k-1}, p_{k-1})\,.
\end{equation*}  Let
the positive constant $\alpha$ with $\beta<\alpha<1$ be determined
later.  We define
\begin{equation*}
\begin{split}
&\mathfrak{R}_{k-1}(u_{k-1},
p_{k-1})\\
&:=(1-\alpha(1-\delta_1))\gamma_1\|\na_{k-1}(u-u_{k-1})\|_{L^2(\Om)}^2
 +\gamma_2(1-\al(1-\delta_2-\delta_3))\|p-p_{k-1}\|_{L^2(\Om)}^2\\[0.5ex]
&\quad+\gamma_3(\beta-\alpha){\eta}^2(u_{k-1}, p_{k-1}, \cT_{k-1})
+\gamma_4(1-\alpha)\sum\limits_{K\in\cT_{k-1}}h_K^{2}\|g\|_{L^2(K)}^2.
\end{split}
\end{equation*}
  Then we perform the  decomposition $\overline{\mathfrak{G}}_{k-1}(u_{k-1}, p_{k-1})=\al\mathfrak{G}_{k-1}(u_{k-1},p_{k-1})+\mathfrak{R}_{k-1}(u_{k-1}, p_{k-1})$ to
   get
\begin{equation*}
\begin{split}
\mathfrak{G}_k(u_k,p_k)
&\leq\al\mathfrak{G}_{k-1}(u_{k-1},p_{k-1})+\mathfrak{R}_{k-1}(u_{k-1}, p_{k-1}).
\end{split}
\end{equation*}
{\bf Step 2~~}
 Now we only need to  show  that it is  possible to choose $\al<1$ such that $\mathfrak{R}_{k-1}(u_{k-1}, p_{k-1})\leq 0$.
 This can be achieved by selecting parameters  $\delta_i\,,i=1,2,3$.  To this end, we recall the reliability of $\eta(u_{k-1}, p_{k-1},  \cT_{k-1})$
  in Lemma \ref{Theorem7.2} with  the multiplication coefficient $C_{Rel}$:
\begin{equation}\label{eq5.8}
\|\na_{k-1}(u-u_{k-1})\|_{L^2(\Om)}^2+\|p-p_{k-1}\|_{L^2(\Om)}^2\leq
C_{Rel}{\eta}^2(u_{k-1}, p_{k-1}, \cT_{k-1})\,.
\end{equation}
Further,  we  take $\delta_1=\delta_2+\delta_3$ with
$0<\delta_1<\min(
\frac{\gamma_3(1-\beta)}{C_{Rel}(\gamma_1+\gamma_2)},1)$.  Then, we
take
$$\alpha:=\frac{(\gamma_1+\gamma_2)C_{Rel}+\gamma_3\beta+\gamma_4}
{(1-\delta_1)(\gamma_1+\gamma_2)C_{Rel}+\gamma_3+\gamma_4}.$$
It is straightforward  to see that $\beta<\al< 1$.  As
\begin{equation}
\sum\limits_{K\in\cT_{k-1}}h_K^{2}\|g\|_{L^2(K)}^2 \leq {\eta}^2(u_{k-1},
p_{k-1}, \cT_{k-1}),
\end{equation}
we obtain
\begin{equation*}
\begin{split}
&\mathfrak{R}_{k-1}(u_{k-1}, p_{k-1})\\
 &\leq
\big((1-\alpha(1-\delta_1))(\gamma_1+\gamma_2)C_{Rel}+\gamma_3(\beta-\alpha)
+\gamma_4(1-\alpha)\big){\eta}^2(u_{k-1}, p_{k-1}, \cT_{k-1})=0\,.
\end{split}
\end{equation*}
 This   proves that
\begin{equation*}
\begin{split}
\mathfrak{G}_k(u_k,p_k)
&\leq\al\mathfrak{G}_{k-1}(u_{k-1},p_{k-1}).
\end{split}
\end{equation*}

{\bf Step 3~~} Finally, we  take $\beta_1:=\gamma_4/\gamma_3$ and
rearrange
$\gamma_2:=\gamma_2(1-\delta_2-\delta_3)/(1-\delta_1)\gamma_1$,
$\gamma_3:=\gamma_3/(1-\delta_1)\gamma_1$, which completes  the proof.
\end{proof}

\section{The discrete reliability}
In this section, we   prove the discrete reliability.  The analysis needs
some prolongation operator from $V_k$ to $V_{k+\ell}$ with some integer $\ell\geq 1$. Some further notations are needed.   Given
 $E\in\cE_{k+\ell}$, the edge patch $\omega_{E, k}$ of $E$ with respect to the
mesh $\cT_k$ is  defined as
\begin{equation}
  \omega_{E, k}:=\{K\in\cT_{k}, E\subset\partial K \text{ or } E \text{ lies in the interior of  }K\}.
 \end{equation}
Let $\xi_{E}=\card( \omega_{E, k})$.  We define the prolongation
interpolation $I_{k+\ell}^{\prime}v_{k}\in V_{k+\ell}$ for any $v_k\in V_k$, as
\begin{equation}\label{eq5.13}
 \int_E I_{k+\ell}^{\prime}v_k ds:=
\frac{1}{\xi_{E}}\sum\limits_{K\in \omega_{E, k}} \int_E (v_k|_{K})
ds\,\text{ for any } E\in\cE_{k+\ell}\,.
\end{equation}
For the interpolation operator $I_{k+\ell}^{\prime}$, we have
\begin{equation}\label{eq5.3new}
I_{k+\ell}^{\prime}v_k=v_k\text{ for any }K\in\cT_k\cap\cT_{k+\ell}\text{ and
}v_k\in V_{k+\ell}.
\end{equation}
As we will see in Remark \ref{remark5.2} below, we cannot  directly use  the prolongation operator $I^{\prime}_{k+\ell}$ in the analysis of the discrete reliability.
An averaging operator is needed.  Denote $\mathcal{N}_k$  as the set of internal vertexes of the mesh $\mathcal{T}_k$, and denote  $S_k\subset H_0^1(\Omega)$ as  the conforming linear element space over $\mathcal{T}_k$.
Given
 $Z\in\mathcal{N}_{k}$, the nodal patch $\omega_{Z, k}$  is defined by
\begin{equation}
  \omega_{Z, k}:=\{K\in\cT_k, Z\in K\}.
 \end{equation}
 Denote $\phi_Z\in S_k$ as  the  canonical basis function associated to $Z$,  which satisfies
$\phi(Z)=1$ and $\phi(Z^{\prime})=0$ for vertex $Z^{\prime}$ of $\cT_k$ other than $Z$.  We define
\begin{equation}
\cE_Z:=\{E\in\cE_k, Z\in\mathcal{N}_k \text{ is one end point of  }E\}.
\end{equation}
The idea of \cite{Brenner2003} leads to the definition of  the following averaging operator $\Pi:  V_k\rightarrow (S_k)^2$:
\begin{equation}
\Pi v_k:=\sum\limits_{Z\in\mathcal{N}_k}v_Z\phi_Z \text{ for any }v_k\in V_k,
\end{equation}
where
\begin{equation}\label{equ5.7}
v_Z=\frac{1}{\xi_Z}\sum\limits_{K\in \omega_{Z, k}}(v_k|_K)(Z)  \text{ with } \xi_{Z}=\card( \omega_{Z, k}).
\end{equation}
Given any $K\in\cT_k$, we have
\begin{equation}\label{eq5.7new}
\begin{split}
&\|\Pi v_k-v_k\|_{L^2(K)}+h_K\|\na(\Pi v_k-v_k)\|_{L^2(K)}\\
&\lesssim h_K^{3/2}(\sum\limits_{T\in\cT_k\& T\cap K\not=\emptyset}\sum\limits_{E\subset\pa T}\|[\na_k v_k\tau_E]\|_{L^2(E)}^2)^{1/2},
\end{split}
\end{equation}
 for any $v_k\in V_k$,  see \cite{Brenner2003} for the  proof.  Define
  \begin{equation*}
 \Omega_{\mathcal{R}}:=\textrm{interior}(\bigcup\{K: K\in\mathcal{T}_k\backslash\mathcal{T}_{k+\ell},\}),
 \end{equation*}
 and
\begin{equation*}
 \Omega_{\mathcal{C}}:=\textrm{interior}(\bigcup\{K: K\in\mathcal{T}_k\cap\mathcal{T}_{k+\ell},\partial K\cap \partial \Omega_{\mathcal{R}}=\emptyset \}).
 \end{equation*}
 The main idea herein is to take the mixture of the prolongation operators $I_{k+\ell}^{\prime}$
  and $\Pi$.  More precisely, we use $\Pi$ in the region $\Omega_{\mathcal{R}}$ where the elements of $\cT_k$ are refined and take $I_{k+\ell}^{\prime}$ in the region $\Omega_{\mathcal{C}}$, and we define some mixture in the layers between them. This leads to the prolongation
   operator $J_{k+\ell}: V_k\rightarrow V_{k+\ell}$ as follows:
   $$
   {J}_{k+\ell}v_k:=\left\{\begin{array}{ll} \Pi v_k &\text{ on }\Omega_{\mathcal{R}},\\ v_k &\text{ on } \Omega_{\mathcal{C}},\\
 v_{k+\ell, tr}&\text{ on }\Omega\backslash(\Omega_{\mathcal{R}}\cup \Omega_{\mathcal{C}}),\end{array}\right.
   $$
 where $v_{k+\ell, tr}$ is defined as
   $$
\int_Ev_{k+\ell, tr}ds:=\left\{\begin{array}{ll}    \int_E\Pi v_kds  &   \text{ for any } E\in \partial (\cT_k\cap\cT_{k+\ell})\\
   \int_E  v_k ds &  \text{ otherwise }.
         \end{array}\right.
   $$
Define
$$
\mathcal{M}_{k, k+\ell}:=\{K\in\mathcal{T}_k, \partial K\cap \overline{\cup(\cT_k\backslash\cT_{k+\ell})}\not=\emptyset \}.
$$

\begin{remark} It follows immediately from regularity of the mesh $\cT_k$ that
$$
\# \mathcal{M}_{k, k+\ell}\leq \kappa \# \cT_k\backslash\cT_{k+\ell}
$$
for a positive constant $\kappa\geq 1$ which is only dependent on the initial mesh $\cT_0$.
\end{remark}
\begin{lemma}\label{Lemma 5.1} For any $v_k\in
V_k$, it holds that
\begin{equation}\label{eq7.72}
\|\nabla_{k+\ell}(J_{k+\ell}v_k-v_k)\|_{L^2(\Om)}^2 \lesssim
\sum\limits_{K\in \mathcal{M}_{k, k+\ell}}\sum\limits_{E\subset \partial K}h_K\|[\nabla_k v_k\tau_E]\|_{L^2(E)}^2\,.
\end{equation}
\end{lemma}
\begin{proof} As $J_{k+\ell} v_k=\Pi v_k$ on $\Omega_{\mathcal{R}}$ and $J_{k+\ell}v_k=v_k$ on $\Omega_{\mathcal{C}}$,
 from \eqref{eq5.3new} and \eqref{eq5.7new}, we
  only need to estimate $\|\na(J_{k+\ell}v_k-v_k)\|_{L^2(K)}=\|\na(v_{k+\ell, tr}-v_k)\|_{L^2(K)}$ for $\cT_{k}\ni K\subset \Omega\backslash(\Omega_{\mathcal{R}}\cup \Omega_{\mathcal{C}} )$.
  Given $E\in\cE_{k}$,  let $\varphi_E$ be the canonical basis function of the nonconforming $P_1$ element on $\cT_{k}$, which satisfies $\int_E\varphi_E ds=|E|$ and $\int_{E^{\prime}}\varphi_E ds=0$
   for any $E^{\prime}\in\cE_{k}$ other than $E$. A direct calculation yields
   \begin{equation*}
   \|\varphi_E\|_{L^2(\Om)}+h_E\|\na_{k}\varphi_E\|_{L^2(\Om)}\lesssim h_E.
   \end{equation*}
   Let $v_{E}^{\prime}:=\int_E v_{k+\ell, tr}|_Kds$ and $v_E:=\int_Ev_k|_Kds$; thus  we have
   \begin{equation}\label{eq5.10new}
   \|\na(v_{k+\ell, tr}-v_k)\|_{L^2(K)}\lesssim \sum\limits_{E\subset\pa K}|v_{E}^{\prime}-v_E|/h_E.
   \end{equation}
   We only need to bound the terms $|v_E^{\prime}-v_E|$ for  $E\subseteq\pa (\cT_k\cap\cT_{k+\ell})$.
    We assume that $Z_1$ and $Z_2$ are two endpoints of $E$.  Then,  the trace of $v_k|_K$ on $E$ can be expressed as
    \begin{equation}
    v_k|_E=(v_k|_K)(Z_1)\phi_{Z_1}+(v_k|_K)(Z_2)\phi_{Z_2}.
    \end{equation}
    Note that
    \begin{equation}
    \Pi v_k|_E=v_{Z_1}\phi_{Z_1}+v_{Z_2}\phi_{Z_2}.
    \end{equation}
    We recall that  $v_{Z_i}$ are  defined in \eqref{equ5.7} and that $\phi_{Z_i}$ are the canonical basis functions
     associated with  vertexes $Z_i$ for  the conforming linear element.
   Therefore
   \begin{equation}\label{eq5.12new}
   \begin{split}
   |v_{E}^{\prime}-v_E|&=|\int_E(\Pi v_k|_E-v_k|_E)ds|\\[0.5ex]
   &=|\int_E ((v_{Z_1}-(v_k|_K)(Z_1))\phi_{Z_1}+(v_{Z_2}-(v_k|_K)(Z_2))\phi_{Z_2})ds|\\
   &\lesssim h_E(\sum\limits_{i=1}^2\sum\limits_{E^{\prime}\in\cE_{Z_i}}h_E^{\prime}\|[\nabla_kv_k\tau_{E^{\prime}}]\|_{L^2(E^{\prime})}^2)^{1/2}.
   \end{split}
   \end{equation}
   By inserting the estimates of  $|v_{E}^{\prime}-v_E|$ in \eqref{eq5.12new} into \eqref{eq5.10new}, we  complete the proof.
  \end{proof}

We define the ratio $\gamma$ as follows:
\begin{equation}\label{ratio}
\gamma:=\max\limits_{K\in\cT_k\backslash\cT_{k+\ell}}\max\limits_{\cT_{k+\ell}\ni T\subset K}\frac{h_K}{h_T}
\end{equation}

\begin{lemma}\label{Lemma7.15} The following discrete reliability holds:
\begin{equation}\label{eq7.60}
\|\na_{k+\ell}( u_{k+\ell}-u_k)\|_{L^2(\Om)}+\|p_{k+\ell}-p_k\|_{L^2(\Om)} \lesssim
{\eta}(u_k, p_k, \mathcal{M}_{k, k+\ell})\,.
\end{equation}
\end{lemma}
\begin{remark}\label{remark5.2}
If we directly take the prolongation operator $I_{k+\ell}^{\prime}$ to analyze this discrete reliability, the constant
for the {\em established } discrete reliability  will depend on the ratio $\gamma$ (see Appendix A for an example).
\end{remark}

\begin{proof} For any $v_{k+\ell}\in V_{k+\ell}$, we have the following
decomposition:
\begin{equation}\label{eq7.61}
\begin{split}
&\mu\|\na_{k+\ell}(u_{k+\ell}-u_k)\|_{L^2(\Om)}^2\\
 &=a_{k+\ell}(u_{k+\ell}-u_k, u_{k+\ell}-v_{k+\ell})+a_{k+\ell}(u_{k+\ell}-u_k, v_{k+\ell}-u_k)\,.
\end{split}
\end{equation}
We will first estimate the first term  on the right--hand side of the
 above  equation.  It follows the discrete problem \eqref{eq7.3} that
\begin{equation}\label{eq7.62}
a_{k+\ell}(u_{k+\ell}-u_k, u_{k+\ell}-v_{k+\ell})
=\res_k(u_{k+\ell}-v_{k+\ell})-b_{k+\ell}(u_{k+\ell}-v_{k+\ell},p_{k+\ell}-p_k)
 \,.
\end{equation}
The first term on the right--hand side of \eqref{eq7.62} can be
bounded as in \eqref{estimateresidual}:
\begin{equation}\label{eq7.68}
\begin{split}
&|\res_k(u_{k+\ell}-v_{k+\ell})| \lesssim
(\sum\limits_{K\in \mathcal{M}_{k, k+\ell}}
 h_K^2\|g\|_{L^2(K)}^2)^{1/2}\|\na_{k+\ell}(u_{k+\ell}-v_{k+\ell})\|_{L^2(\Om)}\,.
 \end{split}
\end{equation}
Now we turn to the second term on the right hand side of
\eqref{eq7.62}.  Thanks to the discrete inf-sup condition, we use
the discrete problem \eqref{eq7.3} to get
\begin{equation}\label{eq7.69}
\begin{split}
\|p_{k+\ell}-p_k\|_{L^2(\Om)}& \lesssim \sup\limits_{0\not=v_{k+\ell}\in
V_{k+\ell}}\frac{b_{k+\ell}(v_{k+\ell},p_{k+\ell}-p_k)}{\|\na_{k+\ell}v_{k+\ell}\|_{L^2(\Om)}}\\[0.5ex]
&\lesssim \sup\limits_{0\not=v_{k+\ell}\in
V_{k+\ell}}\frac{\res_k(v_{k+\ell})}{\|\na_{k+\ell}v_{k+\ell}\|_{L^2(\Om)}}
+\|\na_{k+\ell}(u_{k+\ell}-u_k)\|_{L^2(\Om)}\,.
\end{split}
\end{equation}
 An application of the Cauchy--Schwarz inequality leads to
\begin{equation}\label{eq7.70}
|b_{k+\ell}(u_{k+\ell}-v_{k+\ell},p_{k+\ell}-p_k)|
 \leq \|p_{k+\ell}-p_k\|_{L^2(\Om)}
\|\na_{k+\ell}(u_{k+\ell}-v_{k+\ell})\|_{L^2(\Om)}\,.
\end{equation}
After inserting \eqref{eq7.62}, \eqref{eq7.68}, \eqref{eq7.69}, and
\eqref{eq7.70} into  \eqref{eq7.61},  we use the triangle and
Young inequalities to derive
\begin{equation}\label{eq7.71}
\begin{split}
&\|\na_{k+\ell}(u_{k+\ell}-u_k)\|_{L^2(\Om)}^2
+\|p_{k+\ell}-p_k\|_{L^2(\Om)}^2\\[1.0ex]
& \lesssim \sum\limits_{K\in \mathcal{M}_{k, k+\ell}}
 h_K^2\|g\|_{L^2(K)}^2 +\inf\limits_{v_{k+\ell}\in
 V_{k+\ell}}\|\na_{k+\ell}(u_k-v_{k+\ell})\|_{L^2(\Om)}^2\,.
 \end{split}
\end{equation}
An application of \eqref{eq7.72} bounds the second term on the right--hand side of \eqref{eq7.71}. This completes the proof.
\end{proof}

With $\gamma_1$ from Theorem \ref{Theorem7.8}, we define the
following  energy norm:
\begin{equation}\label{Norm3b}
\Vvert v,q\Vvert^2:=\|\na
v\|_{L^2(\Om)}^2+\gamma_1\|q\|_{L^2(\Om)}^2, \text{ for any }
(v,q)\in V\x Q.
\end{equation}
We denote  its piecewise version by $\Vvert\cdot\Vvert_{k+\ell}$.

The following lemma gives  links between the error reduction to the
bulk criterion.
\begin{lemma}\label{Theorem7.16}
Let $\cT_{k+\ell}$ be the refinement of $\cT_k$ with the following
reduction:
\begin{equation}\label{eq7.73}
\begin{split}
\Vvert  {u}-u_{k+\ell},  {p}-p_{k+\ell}\Vvert^2_{k+\ell}+\gamma_2\osc^2(g,\cT_{k+\ell})\\
 \leq
\al^{\prime}\big( \Vvert
 {u}-u_k,  {p}-p_k\Vvert^2_k+\gamma_2\osc^2(g, \cT_k)\big),
\end{split}
\end{equation}
with $0<\al^{\prime}<1$ and  the positive constant  $\gamma_2$ from
Theorem \ref{Theorem7.8}. There exists $0<\theta_{\ast}<1$ with
\begin{equation}\label{eq7.74}
\theta_{\ast}\eta^2(u_k, p_k,\cT_k)\leq\eta^2(u_k,p_k,
{\mathcal{M}_{k, k+\ell}}).
\end{equation}
\end{lemma}
\begin{proof} It follows \eqref{eq7.73} and the definitions of the norms $\Vvert\cdot\Vvert_k$ and $\Vvert\cdot\Vvert_{k+\ell}$ that
\begin{equation*}
\begin{split}
& (1-\al^{\prime})(\Vvert  {u}-u_k, {p}-p_k\Vvert^2_k+\gamma_2\osc^2(g,\cT_k))\\[0.5ex]
&\leq\Vvert
 {u}-u_k, {p}-p_k\Vvert^2_k+\gamma_2\osc^2(g,\cT_k)-\Vvert
 {u}-u_{k+\ell}, {p}-p_{k+\ell}\Vvert^2_{k+\ell}-\gamma_2\osc^2(g,\cT_{k+\ell})\\[0.5ex]
&=\|\na_{k+\ell}( u_k-u_{k+\ell})\|_{L^2(\Om)}^2
+\gamma_1\|p_k-p_{k+\ell}\|_{L^2(\Om)}^2+\frac{2}{\mu}a_{k+\ell}( {u}-u_{k+\ell},
u_{k+\ell}-u_k)\\[0.5ex]
&\quad+2\gamma_1( {p}-p_{k+\ell},
p_{k+\ell}-p_k)_{L^2(\Om)}+\gamma_2\osc^2(g,\cT_k)-\gamma_2\osc^2(g,\cT_{k+\ell})\\[0.5ex]
&=I_1+I_2+I_3+I_4+I_5.
 \end{split}
\end{equation*}
The first two terms, $I_1$ and $I_2$,  are estimated  by the discrete
reliability in Lemma \ref{Lemma7.15},
\begin{equation}\label{eq5.14b}
\begin{split}
\Vvert u_{k+\ell}-u_k\Vvert_{k+\ell}^2+\gamma_1\|p_k-p_{k+\ell}\|_{L^2(\Om)}^2\leq
C_{Drel} \eta^2(u_k,p_k, {\mathcal{M}_{k, k+\ell} }),
\end{split}
\end{equation}
where the coefficient $C_{Drel}$ is from Lemma \ref{Lemma7.15}.   The third term $I_3$ can be estimated by the quasi--orthogonality
   of  the velocity in Lemma \ref{Lemma7.5}.  In fact, let the multiplication constant  therein   be the coefficient $C_{QOV}$, so that we have
\begin{equation}\label{equ5.26}
\begin{split}
 &|\frac{2}{\mu}a_{k+\ell}(u-u_{k+\ell},u_{k+\ell}-u_k)|\\
 &\leq 2C_{QOV}\|\na_{k+\ell}(u-u_{k+\ell})\|_{L^2(\Om)} \big(\sum\limits_{K\in
  \mathcal{M}_{k, k+\ell} }h_K^{2}\|g\|_{L^2(K)}^2\big)^{1/2}\\[0.5ex]
  &\leq
  \frac{1-\al^{\prime}}{2}\|\na_{k+\ell}(u-u_{k+\ell})\|_{L^2(\Om)}^2+\frac{2(C_{QOV})^2}{1-\al^{\prime}}\sum\limits_{K\in
  \mathcal{M}_{k, k+\ell} }h_K^{2}\|g\|_{L^2(K)}^2.
\end{split}
\end{equation}
Next,  we use the quasi--orthogonality of  the pressure in Lemma \ref{Lemma7.6} to analyze the fourth term, $I_4$.
Denote the constant of Lemma \ref{Lemma7.6} by $C_{QOP}$, and we obtain
\begin{equation*}
\begin{split}
 &|2\gamma_1(p-p_{k+\ell}, p_{k+\ell}-p_k)_{L^2(\Om)}|\\[0.5ex]
 & \leq 2\gamma_1C_{QOP}
 \bigg(\big(\sum\limits_{K\in \mathcal{M}_{k, k+\ell}}
h_K^2\|g\|_{L^2(K)}^2\big)^{1/2}
 +\|\na_{k+\ell}(u_{k+\ell}-u_k)\|_{L^2(\Om)}\bigg)\|p-p_{k+\ell}\|_{L^2(\Om)}\\[0.5ex]
 & \leq\frac{2\gamma_1(C_{QOP})^2}{1-\al^{\prime}}\bigg(\big(\sum\limits_{K\in \mathcal{M}_{k, k+\ell}}
h_K^2\|g\|_{L^2(K)}^2\big)^{1/2}
 +\|\na_{k+\ell}(u_{k+\ell}-u_k)\|_{L^2(\Om)}\bigg)^2\\[0.5ex]
 &\quad + \frac{1-\al^{\prime}}{2}\gamma_1\|p-p_{k+\ell}\|_{L^2(\Om)}^2.
\end{split}
 \end{equation*}
 Hence it follows from \eqref{eq5.14b} that
\begin{equation}
\begin{split}
&|2\gamma_1(p-p_{k+\ell}, p_{k+\ell}-p_k)_{L^2(\Om)}|\\[0.3ex]
 &\leq \frac{1-\al^{\prime}}{2}\gamma_1\|p-p_{k+\ell}\|_{L^2(\Om)}^2+
 \frac{2\gamma_1(C_{QOP})^2(1+C_{Drel}^{1/2})^2}{1-\al^{\prime}}\eta^2(u_k,p_k, \mathcal{M}_{k, k+\ell}).
 \end{split}
\end{equation}
A direct calculation leads to
 \begin{equation}\label{eq5.27}
 \gamma_2|\osc^2(f,\cT_k)-\osc^2(f, \cT_{k+\ell})|\leq \gamma_2\eta^2(u_k,
 p_k, \mathcal{M}_{k, k+\ell})\,,
 \end{equation}
 we combine \eqref{eq5.14b}---\eqref{eq5.27},  and \eqref{eq7.73} with the efficiency of
 the estimator, which   proves the desired result by
 the parameter
 $$\theta_{\ast}=\frac{(1-\al^{\prime})^2C_{Eff}}{2(2(C_{QOV})^2+2\gamma_1(C_{QOP})^2(1+C_{Drel}^{1/2})^2
 +(1-\al^{\prime})(C_{Drel}+\gamma_2))}\,,
 $$
 with  the efficiency constant $C_{Eff}$ of the estimator $\eta(u_k,
 p_k, \cT_k)$ from Lemma \ref{Theorem7.2}.
\end{proof}

\section{The optimality of the ANFEM}\label{Sec8}
 In this section, we   address the  optimality  of
the adaptive nonconforming linear element method  under
consideration.  We need to control  the consistency error $\consis(\sigma,\cT)$  defined by
\begin{equation}\label{eq7.54}
\consis(\sigma,\cT)=\sup\limits_{v_{\cT}\in V_{\cT}}\frac{(g,
v_{\cT})_{L^2(\Om)}-(\sigma, \na_{\cT}
v_{\cT})_{L^2(\Om)}}{\|\na_{\cT}v_{\cT}\|_{L^2(\Om)}}\quad \text{ with }
\sigma=\mu\na u+p\id\,,
\end{equation}
where $\cT$ is some refinement of the initial mesh $\cT_0$ by the newest vertex bisection.
The following conforming finite element space is needed:
\begin{equation}
P_3(\cT):=\{v\in (H_0^1(\Om))^2, v|_K\in (P_3(K))^2, \text{ for any }K\in\cT\}.
\end{equation}
Then,  there exists an interpolation operator $\Pi_{\cT}: V_{\cT}\rightarrow P_3(\cT)$  with
 the following properties \cite[Lemma A.3]{GiraultRaviart1986}:
\begin{equation}\label{PropertyIn}
\begin{split}
&\int_E (v_{\cT}-\Pi_{\cT}v_{\cT})\cdot c_Eds=0 \text{ for any }c_E\in (P_1(E))^2,\\
&\int_K (v_{\cT}-\Pi_{\cT}v_{\cT}) dx=0,
\end{split}
\end{equation}
for any edge $E$ and element $K$ of $\cT$. In addition, we have
\begin{equation}\label{PropertyIn2}
\|v_{\cT}-\Pi_{\cT}v_{\cT}\|_{L^2(K)}+h_K\|\na\Pi_{\cT}v_{\cT}\|_{L^2(K)}\lesssim h_K \|\na_{\cT}v_{\cT}\|_{L^2(\omega_K)}.
\end{equation}
For any $s_{\cT}\in V_{\cT}$ and $q_{\cT}\in Q_{\cT}$, we define $\sigma_{\cT}=\mu s_{\cT}+q_{\cT}$.  The idea of \cite[Lemma 2.1]{Gudi2010}
 leads to the following decomposition:
\begin{equation}
\begin{split}
&(g,
v_{\cT})_{L^2(\Om)}-(\sigma, \na_{\cT}
v_{\cT})_{L^2(\Om)}\\
&=(g,
v_{\cT}-\Pi_{\cT}v_{\cT})_{L^2(\Om)}-(\sigma-\sigma_{\cT}, \na_{\cT}
(v_{\cT}-\Pi_{\cT}v_{\cT}))_{L^2(\Om)}\\
&\quad + (\sigma_{\cT}, \na_{\cT}
(v_{\cT}-\Pi_{\cT}v_{\cT}))_{L^2(\Om)}
\end{split}
\end{equation}
for any $v_{\cT}\in V_{\cT}$.  By the properties \eqref{PropertyIn} and \eqref{PropertyIn2}, we obtain
\begin{equation}\label{consis}
\consis(\sigma,\cT)\lesssim \inf\limits_{(v_{\cT},q_{\cT})\in
V_{\cT}\x Q_{\cT}}\Vvert u-v_{\cT},
p-q_{\cT}\Vvert_{\cT}+\osc(g,
\cT).
\end{equation}
This implies that the nonlinear approximate class used in \cite{HuXu2007} is {\em equivalent } to the standard nonlinear
 approximate class \cite{CasconKreuzerNochettoSiebert07}.
Hence, we can  introduce the following semi-norm:
\begin{equation}\label{eq7.53}
\mathfrak{E}^2(N;u,p,
g):=\inf\limits_{\cT\in\mathbb{T}_N}\big(\inf\limits_{(v_{\cT},q_{\cT})\in
V_{\cT}\x Q_{\cT}}\Vvert u-v_{\cT},
p-q_{\cT}\Vvert^2_{\cT}+\gamma_2\osc^2(g,
\cT)\big)\,.
\end{equation}
Then the nonlinear approximate class $\mathbb{A}_s$ can be defined
 by
\begin{equation}\label{eq7.55}
\mathbb{A}_s:=\{(u,p, g), |u,p,
g|_s:=\sup\limits_{N>N_0}N^{s}\mathfrak{E}(N;u,p, g)<+\infty\}.
\end{equation}

We must stress that this is the first time  the standard nonlinear approximate class \cite{CasconKreuzerNochettoSiebert07}
has been used to analyze  the adaptive nonconforming finite element method.  In the relevant literature,
the discrete solution of the discrete problem has been  used to define the nonlinear approximate class \cite{BeckerMaoShi2010,BeckerMao2011,MaoZhaoShi2010,Rabaus2010}.  Let $(u_{\cT}, p_{\cT})$ be the approximation solution
of \eqref{eq7.3} on the mesh $\cT$. It follows from the Strang Lemma \cite{CiaBook}
$$
\Vvert u-u_{\cT}, p-p_{\cT}\Vvert_{\cT}\lesssim \inf\limits_{(v_{\cT},q_{\cT})\in
V_{\cT}\x Q_{\cT}}\Vvert u-v_{\cT},
p-q_{\cT}\Vvert_{\cT}+\consis(\sigma,\cT),
$$
and the following fact
$$
\inf\limits_{(v_{\cT},q_{\cT})\in
V_{\cT}\x Q_{\cT}}\Vvert u-v_{\cT},
p-q_{\cT}\Vvert_{\cT}+\consis(\sigma,\cT)\lesssim  \Vvert u-u_{\cT}, p-p_{\cT}\Vvert_{\cT},
$$
that the nonlinear approximate class of \cite{BeckerMao2011} is equivalent to $\mathbb{A}_s$ of \eqref{eq7.55}.
A similar method herein proves that the nonlinear approximate class of  \cite{BeckerMaoShi2010,MaoZhaoShi2010,Rabaus2010}
is equivalent to the standard nonlinear approximate class \cite{CasconKreuzerNochettoSiebert07}.

\begin{remark} After we submitted the revised version to the journal,
 we learnt about that a different argument of \cite{CarstensenPeterseimSchedensack2012} shows that the nonlinear approximate class of  \cite{BeckerMaoShi2010,MaoZhaoShi2010,Rabaus2010}
  is  equivalent to the standard nonlinear approximate class \cite{CasconKreuzerNochettoSiebert07}.
\end{remark}

Thanks to \eqref{consis}, we have
\begin{equation}\label{eq7.79}
\begin{split}
&\Vvert  {u}-u_{k-1},  {p}-p_{k-1}\Vvert^2_{k-1}\\
& \lesssim \inf\limits_{(v_{k-1},
q_{k-1})\in V_{k-1}\x Q_{k-1}}\Vvert  {u}-v_{k-1},
 {p}-q_{k-1}\Vvert^2_{k-1}+\osc^2(g, \cT_{k-1})\,.
\end{split}
\end{equation}
A straightforward investigation shows that if $\cT_k$ is any
refinement of $\cT_{k-1}$,  then
it  holds that
\begin{equation}\label{eq7.75}
\begin{split}
&\inf\limits_{(v_k, q_k)\in V_k\x Q_k}\Vvert
 {u}-v_k,
 {p}-q_k\Vvert^2_k+\gamma_2\osc^2(g, \cT_k)\\[0.5ex]
&\qquad\leq C_{3}\big(\inf\limits_{(v_{k-1}, q_{k-1})\in V_{k-1}\x Q_{k-1}}\Vvert
 {u}-v_{k-1},  {p}-q_{k-1}\Vvert^2_{k-1}+\gamma_2\osc^2(g, \cT_{k-1})\big)\,.
\end{split}
\end{equation}

With these preparations,  following
\cite{HSX10}, we have the following optimality:

\begin{theorem}\label{Theorem7.20}
 Let
$( {u}, {p})$ be the solution of Problem \eqref{eq7.1}, and let
$(\cT_k, V_k\x Q_k, (u_k,p_k))$ be the sequence of meshes, finite
element spaces, and discrete solutions produced by the adaptive
finite element methods. If $( {u}, {p}, g)\in \mathbb{A}_s$ with
$$\theta\leq\frac{C_{Eff}}{2(2(C_{QOV})^2+2\gamma_1(C_{QOP})^2(1+C_{Drel}^{1/2})^2
 +C_{Drel}+\gamma_2)}\,.
 $$
 Then,  it holds that
\begin{equation}\label{eq7.87}
\Vvert u-u_N, p-p_N\Vvert_N^2+\gamma_2\osc^2(g, \cT_N)\lesssim |
{u},
 {p}, g|_s^2(\#\cT_N-\#\cT_0)^{-2s}\,.
\end{equation}
\end{theorem}

\appendix
\section{A counter  example}
We present an example in this appendix to show that if the prolongation operator  $I_{h}^{\prime}$
  defined  by \eqref{eq5.13} is directly used to analyze the discrete reliability of the estimator,  the  constant for the {\em established }
  discrete reliability  could depend on some key mesh
refinement ratio
$$\gamma:=\max\limits_{K\in\cT_H\backslash\cT_h}\max\limits_{\cT_h\ni
  T\subset K}\frac{h_K}{h_T},
  $$
  where $\cT_H$ is some regular  triangulation of $\Om$ into triangles and $\cT_h$ is some refinement of
 $\cT_H$.  To this end, we first give an example to demonstrate  that  there are generally  no positive constants $C$  independent of $\gamma$ such that
the following estimate holds true:
\begin{equation}\label{counterexample}
\sum\limits_{E\in\cE_h\backslash\cE_H}\int_E[u_H]\{\frac{\pa v_h}{\pa\nu_E}\}ds
\leq C(\sum\limits_{E\in\cE_H\backslash\cE_h}h_E^{-1}\|[u_H]\|_{L^2(E)}^2)^{1/2}\|\na_hv_h\|_{L^2(\Om)},
\end{equation}
where $u_H\in V_H$ is the  finite element solution of the velocity on the mesh $\cT_H$ and $v_h$ is some element of $V_h$ over the nested fine mesh $\cT_h$.
As usual,  $\cE_H$ (resp. $\cE_h$) is the set of the edges of
$\cT_H$ (resp. $\cT_h$). Denote $V_H$ (resp. $V_h$) as  the nonconforming
linear element space with respect to $\cT_H$
(resp. $\cT_h$). Denote $[\cdot]$ as the jump of some function across the
edge $E$ and $\{\cdot\}$ as the average of some function across the
edge $E$.  In addition, denote $\nu_E$ as the unit normal vector to $E$ with
the length $h_E$.

 In the  following,  an example is given to show that  $u_H\in V_H$ and $v_h\in V_h$ exist such that
  the above constant $C$ depends on the ratio $\gamma$.
  For simplicity,  let $\cT_H$ consist of two triangles $\triangle ABC$ and $\triangle ACD$ as in  Figure \ref{figure}.
Let $\cT_h$ be a uniform triangulation of $\Om$ into $2\times N^2$ triangles, cf. Figure \ref{figure} for the case $N=5$.  We stress that the idea and result  can be easily extended to the mesh with  the newest vertex bisection. For the sake of simplicity,
let $N=2k+1$ with some nonnegative integer $k$. Let $Z_i$, $i=-k,\cdots, k$, be the nodes of $\cT_h$ whose coordinates are $(\frac{1}{N},\frac{2i}{N})$. Let $\phi_{Z_i}$ be the nodal basis function of the conforming linear element space defined over $\cT_h$ such that
  $\phi_{Z_i}(Z_i)=1$ and $\phi_{Z_i}(Z)=0$ for any node $Z$ other than $Z_i$.  We choose $u_H\in V_H$ such that the jump is $[u_H]=y$ over the edge $AC$. We choose $v_h$ as follows:
  \begin{equation}
  v_h:=\sum\limits_{i=-k}^{k}sign(i)\phi_{Z_i}\text{ with } sign(i):=\left\{\begin{split}&1\text{ if }i>0,\\ & 0 \text{ if }i=0,\\&-1\text{ if }i<0.\end{split}\right.
  \end{equation}
\begin{figure}[h]
\setlength{\unitlength}{0.40cm}
\begin{picture}(20,20)
\thicklines
\thicklines
\put (10,0){\line(0,1){20}}
\put(10,0){\line(1,1){10}}
\put(10,0){\line(-1,1){10}}
\put(10,20){\line(1,-1){10}}
\put(10,20){\line(-1,-1){10}}
\thinlines
\multiput (10,0)(2,2){5}{\line(-1,1){10}}
\multiput (10,0)(-2,2){5}{\line(1,1){10}}
\put(12,2){\line(0,1){16}}
\put(14,4){\line(0,1){12}}
\put(16,6){\line(0,1){8}}
\put(18,8){\line(0,1){4}}
\put(8,2){\line(0,1){16}}
\put(6,4){\line(0,1){12}}
\put(4,6){\line(0,1){8}}
\put(2,8){\line(0,1){4}}
\put(10,-1){A(0,-1)}
\put(21,10){B(1,0)}
\put(10,20){C(0,1)}
\put(-3,10){D(-1,0)}
\end{picture}
\caption{The meshes $\cT_H$ and $\cT_h$.}\label{figure}
\end{figure}
Note that $\{\frac{\pa\phi_{Z_i}}{\pa\nu_E}\}=N/2$ over the edge $AC$ for $i=-k, \cdots, k$.  A direct calculation gives
\begin{equation}
\int_{AC}[u_H]\{\frac{\pa v_h}{\pa \nu_E}\}ds=N/2-\frac{1}{2N}.
\end{equation}
On the other hand, a direct calculation leads to
\begin{equation}
\|\na_hv_h\|^2_{L^2(\Om)}\leq 4N.
\end{equation}
This indicates that the constant $C$ in \eqref{counterexample} should be  $\mathcal{O}(\sqrt{N})$, which depends on the ratio $\gamma=\mathcal{O}(N) $ for this example.

 For the analysis of the discrete reliability,  a direct application of the prolongation operator $I_{h}^{\prime}$ as defined in \eqref{eq5.13} will
    lead to a similar  estimate like \eqref{counterexample}, and , as a result, the constant for the {\em established } discrete reliability based on  such an estimate  will depend on the ratio $\gamma$. Note that in the analysis of
optimality of the adaptive method it is possible to know that $\cT_h$ is
some refinement of $\cT_H$ only by the newest vertex bisection \cite{BinevDahmenDeVore04,Stevenson06,CasconKreuzerNochettoSiebert07}. Note, too, that
 there is no guarantee that  $\gamma$  is bounded. Therefore, the proof of the discrete reliability
 based on  the  prolongation operator $I_h^{\prime}$ as presented in  \cite{BeckerMao2011,HuXu2007} may not lead to a uniform estimate as claimed.

\end{document}